\documentclass[12pt,twocolumn]{amsart}
\usepackage[margin=.95in]{geometry}
\usepackage[tt=false]{libertine}
\usepackage{mathabx}
\usepackage[T1]{fontenc}
\usepackage{microtype}
\usepackage[normalem]{ulem}
\usepackage[dvipsnames]{xcolor}
\usepackage[
  colorlinks=true,
  linkcolor=black,
  citecolor=black,
  urlcolor=black
]{hyperref}
\usepackage{enumerate}

\usepackage{tikz}
\usetikzlibrary{calc} 
\usepackage{pgfplots}
\pgfplotsset{compat=1.18}

\usepackage[backend=biber,style=numeric]{biblatex}
\newlength{\bibcontindent}
\defbibenvironment{bibliography}
  {\list
     {\printtext[labelnumberwidth]{%
        \printfield{labelprefix}%
        \printfield{labelnumber}}}
     {\setlength{\labelwidth}{\labelnumberwidth}%
      \setlength{\leftmargin}{\bibcontindent}%
      \setlength{\itemindent}{\labelwidth}%
      \addtolength{\itemindent}{\labelsep}%
      \addtolength{\itemindent}{-\bibcontindent}%
      \setlength{\itemsep}{\bibitemsep}%
      \setlength{\parsep}{\bibparsep}}%
      }
  {\endlist}
  {\item}

\newcommand{\R}{\mathbb{R}}
\newcommand{\de}{\delta}

\theoremstyle{definition}

\title[What is... An Exceptional Set Estimate?]{What is... \\ An Exceptional Set Estimate?}
\author{Paige Bright}
\author{Krystal Taylor}

\begin{document}

\maketitle

\section*{Shadows are Typically Large}

\noindent Suppose you are standing outside on a sunny day from dawn til dusk, observing your shadow as it changes. Most of the time, you will see a large portion of your shadow; but at noon, when the sun is right above you, your shadow will fall right beneath your feet. In this way, at almost every time of day, your shadow is relatively large, and, more quantitatively, it is very \textit{exceptional} for your shadow to be small. In projection theory, a field at the intersection of discrete geometry and geometric measure theory, one seeks to understand how general this phenomenon is in a mathematically precise way.

The concept of shadow here is made precise by orthogonal projections. Given $\theta \in \mathbb{S}^1 \subset \R^2$, we let $P_\theta : \R^2 \to \R$ denote the orthogonal projection onto the line through the origin in direction $\theta$. In other words, $P_\theta(x) = \theta \cdot x$ where $\cdot$ is the Euclidean dot product. Given $X\subset \R^2$ Borel, how large could the projection, $P_\theta(X)$, be? Answers to this question depend on how large the underlying set $X$ is itself. By Fubini's theorem, if $X$ has positive planar Lebesgue measure, then $P_\theta(X)$ must have positive linear (1-dimensional) Lebesgue measure for all $\theta$. As such, we consider sets which have null Lebesgue measure but that are large in another quantitative sense, namely Hausdorff dimension. 

Before continuing, we briefly note that the Hausdorff dimension of a set $X$ in a metric space, denoted $\dim X$ throughout, can take on any real nonnegative value bounded by the dimension of the ambient space. Hausdorff dimension is appropriate for quantifying the size of Lebesgue measure zero sets. For instance, while the standard middle thirds Cantor set has Lebesgue measure zero, it has Hausdorff dimension $\frac{\log 2}{\log 3}$. Heuristically, $\dim X =\alpha$ if the minimum number of $\de$-balls needed to cover $X$, $|X|_\de$, is $\sim \de^{-\alpha}$ for sufficiently small $\de$. This heuristic would be made precise if we were instead using \textit{box dimension}, but this heuristic is sufficient for our purposes. For more background on Hausdorff dimension, see e.g. \cite{Mattila95}.

Using the fact that orthogonal projection is a Lipschitz map onto a 1-dimensional subspace, it is immediate to see that 
\[
\dim P_\theta(X) \leq \min\{\dim X, 1\} \,\, \text{for all $\theta \in \mathbb{S}^1$}.
\]
Intuitively, all this is saying is that the size of the shadow ($\dim P_\theta(X)$) is at most the size of the object casting the shadow ($\dim X$), and at most the size of what the shadow is being cast upon ($\dim \R = 1$). In 1954, Marstrand \cite{Marstrand54} showed that in fact,
\[
\dim P_\theta(X) = \min\{\dim X,1 \} \,\, \text{for a.e. $\theta \in \mathbb{S}^1.$}
\]
(Here, almost every is with respect to the surface measure on the circle.) In other words, \textit{almost all of the time}, the orthogonal projection of $X$ is \textit{as large as it could be}, matching our intuition from our nice sunny day outside. 

However, projection theorists sought to probe this question even further by estimating the size of the set of $\theta$ where the above equality fails. These are known as \textit{exceptional set estimates}.

\section*{A Discrete Model and Sharpness}

It is natural to first approach these estimates by studying a discrete model where $X \subset \R^2$ is finite. 
To begin, 
it is immediate to show that for all but finitely many $\theta$, 
\begin{equation} \label{eqn:finiteMarstrand}
|P_\theta(X)| = |X|.
\end{equation}
Indeed, in order for $P\theta(X)$ to be smaller, there must be two distinct points $x,y \in X$ such that $\theta \cdot (x-y) =0$, and $X$ can only determine at most $\binom{|X|}{2}$ such directions. Thus, let us bound the number of $\theta$ for which \eqref{eqn:finiteMarstrand} fails. In particular, given $s \leq |X|$, we seek to bound the cardinality of 
\[
E_s(X) := \{\theta \in \mathbb{S}^1_+ : |P_\theta(X)| < s\}.
\]
We consider $\mathbb{S}^1_+$, i.e. angles between $[0,\pi)$, so as to avoid double-counting.
As a starting point, let's consider two model examples of sets $X\subset \R^2$. If $X$ is a collection of $N$ points on the $y$-axis, then $E_s(X)$ is simply the direction $\theta = (1,0)$ for all $s\in (1,N]$. This matches our sunny day picture: The only time in which $X$ can cast an exceptionally small shadow is at ``noon.'' If $X$ is a \textit{generic set} of $N$ points, which is to say no line contains more than 2 points of $X$, then 
\[
E_{s}(X) = \emptyset \quad \text{for } s\leq N/2.
\]
To see this, notice that for all $\theta \in E_s(X)$, by the pigeonhole principle, one of the fibers of $P_\theta$ must contain at least three points of $X$ which is a contradiction. 
For the interested reader, we leave as an exercise to consider $E_{\lceil N/2\rceil+1}(X)$ when $X$ is a regular $N$-gon, in comparison to the set of directions between pairs of points in $X$.

Let us now prove bounds on $E_s(X)$ for general $X\subset \R^2$. One estimate we can note is that 
\begin{equation}\label{eqn:product}
|E_{|X|^{1/2}}(X)| \leq 1.
\end{equation}
Indeed, if $|E_{|X|^{1/2}} (X)| \geq 2$, there exist two nonantipodal directions $\theta,\theta'$ such that we can cover $X$ in $< |X|^{1/2}$ many parallel lines in directions $\theta^\perp$ and $(\theta')^\perp$ respectively. We denote these lines/fibers $\mathcal L_\theta$ and $\mathcal L_{\theta'}$. Since both sets of lines cover $X$, it follows that 
\[
|X| \leq |\mathcal L_{\theta}||\mathcal L_{\theta'}| < |X|^{1/2} \cdot |X|^{1/2} = |X|,
\]
which is a contradiction.

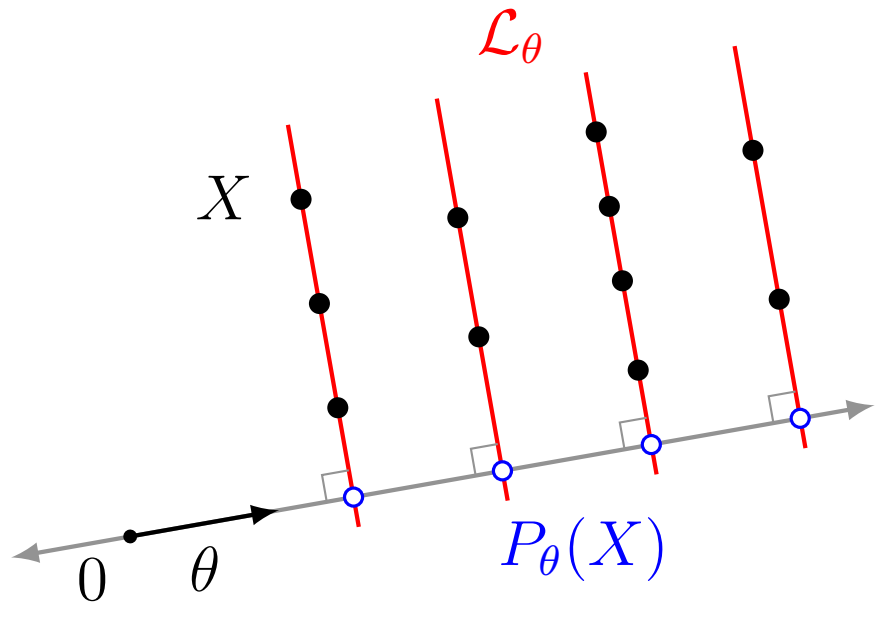
\begin{figure} 
    \centering
    \begin{tikzpicture}[
  fib/.style={red, thick},
  pt/.style={circle, fill=black, inner sep=1.4pt},
  proj/.style={circle, draw=blue, semithick, fill=white, inner sep=1.2pt}
]
\def\H{2.5}  
\begin{scope}[rotate=10]
  \draw[gray!85, thick, latex-latex] (-0.8,0) -- (5,0);
  \foreach \c in {1.5,2.5,3.5,4.5} {
    \draw[fib] (\c,-0.2) -- (\c,\H);
    \draw[thin, gray!85] (\c-0.18,0) -- (\c-0.18,0.18) -- (\c,0.18);
  }
  \foreach \c in {1.5,2.5,3.5,4.5}
    \node[proj] at (\c,0) {};
  \node[blue, below] at (3.0,-0.25) {$P_\theta(X)$};
  \node[red, above] at (3,\H) {$\mathcal{L}_\theta$};
  \foreach \x/\t in {1.5/0.24, 1.5/0.52, 1.5/0.80,
                     2.5/0.36, 2.5/0.68,
                     3.5/0.20, 3.5/0.44, 3.5/0.64, 3.5/0.84,
                     4.5/0.32, 4.5/0.72}
    \node[pt] at (\x,\t*\H) {};
  \node at (1.0,\H-.4) {$X$};
  \fill (0,0) circle (1.3pt) node[below left] {$0$};
  \draw[-latex, thick] (0,0) -- (1,0) node[midway, below] {$\theta$};
\end{scope}
\end{tikzpicture}
    \caption{Covering $X$ by fibers $\mathcal L_\theta$.}
    \label{figure1}
\end{figure}

To obtain stronger estimates, we consider an incidence problem between points and lines. To do so, first note that by definition, for all $\theta \in E_s(X)$, it takes at most $s$-many lines in direction $\theta^\perp$ to cover $X$. Call this set of fibers $\mathcal L_\theta$ (see Figure \ref{figure1}), and let 
\[
\mathcal L := \bigcup_{\theta \in E_s(X)} \mathcal L_\theta.
\]
Then, we consider the set of \textit{incidences} between $X$ and $\mathcal L$; that is,
\[
\mathcal I(X,\mathcal L) := \{(x,\ell) \in X \times \mathcal L : x \text{ lies in } \ell\}.
\]
\vspace{-.6cm}

Notice that since $\mathcal L_\theta$ covers $X$ and $\mathcal L_\theta \cap \mathcal L_{\theta'} = \emptyset$ for all $\theta \neq \theta'$, it follows that 
\[
|X||E_s(X)| = \sum_{x\in X} |\mathcal I(x, \mathcal L)| = |\mathcal I(X,\mathcal L)|.
\]
Indeed, notice that for every $x\in X$ and every $\theta \in E_s(X)$, there is a distinct fiber in $\mathcal L_\theta$ containing $x$. Hence, if we have an upper bound on the right hand side, we can hope to conclude estimates on $|E_s(X)|$. For instance, using that two points determine one line through both, one can use Cauchy--Schwarz to verify that 
\[
|\mathcal I(X,\mathcal L)| \leq |X||\mathcal L|^{1/2} + |\mathcal L|.
\]
Since $|\mathcal L_\theta| \leq s$, one has $|\mathcal L| \leq s|E_s(X)|$. Hence, if $s<|X|/2$ (and thus the first term dominates) one can conclude
\begin{equation} \label{eqn:CSineq}
|E_s(X)| \lesssim s.
\end{equation}
While \eqref{eqn:CSineq} is not sharp, it does mark a significant improvement on the trivial estimate $$|E_s(X)| \leq \binom{|X|}{2}\sim |X|^2.$$

The sharp estimate was given by Szemer\'edi and Trotter \cite{SzemerediTrotter}, who proved
\[
|\mathcal I(X,\mathcal L)|\leq C (|X|^{2/3} |\mathcal L|^{2/3} + |X| + |\mathcal L|)
\]
for some universal constant $C$. Hence, assuming $s < |X|/3C,$ it follows that 
\begin{equation} \label{eqn:STSharp}
|E_s(X)| \lesssim \max\{s^2 |X|^{-1}, 1\},
\end{equation}
which is sharp. We now turn to the continuous analogues of these problems in the framework of geometric measure theory.

\section*{The Problem and Its History}

Let us begin by stating what the continuum analogue of $E_s(X)$ is: Given $X\subset \R^2$ and some $0<s < \min\{\dim X, 1\}$, let 
\[
E_s(X) = \{\theta \in \mathbb{S}^1_+ : \dim P_\theta(X) \leq s\}.
\]
The aforementioned result of Marstrand shows $E_s(X)$ has null measure for all $s.$ Hence, it is natural to ask for upper bounds on the dimension of $E_s(X)$.


Each of the discrete bounds we derived has an interesting continuum counterpart, with the analogue of \eqref{eqn:STSharp} being conjectured by D. Oberlin in 2008 \cite{Oberlin2s-aConjecture}. To go from the discrete to the continuum, let us go back to our heuristic that $|X|_\de \sim \de^{-\dim X}$ for all small $\de$. Based on the discrete setting, we might expect that $\de$-balls should be no different from points for small enough $\de$. As such, by \eqref{eqn:STSharp} we might conjecture
\begin{align*}
|\{\theta : |P_\theta(X)|_\de \lesssim \de^{-s}\}|_\de &\lesssim \max\{(\de^{-s})^2 |X|_\de^{-1}, \de^{0}\} \\
&\sim \de^{- \max\{2s-\dim X, \,0\}}.
\end{align*}
This thus naturally leads to the following conjecture of Oberlin: For $X\subset \R^2$ Borel, 
\begin{equation}\label{eqn:OSRW}
\dim E_s(X) \leq \max\{2s - \dim X,0\}.
\end{equation}
He obtained the statement for $\dim X \leq 1$ and the special case $s = \frac{\dim X}{2} - \epsilon$ for all $\epsilon>0$ using restricted Radon transforms, which is what motivated this conjecture. The case $ s= \frac{\dim X}{2}$ for the entire range of $\dim X \in [0,2]$ was obtained by Bourgain in 2010 \cite{BourgainProjectionTheorem} and can be thought of as an analogue of \eqref{eqn:product}.

The analogue of \eqref{eqn:CSineq} was proven by Kaufman in 1968 \cite{Kaufman68} who proved 
\begin{equation} \label{eqn:Kaufman}
\dim E_s(X) \leq s.
\end{equation}
Notice however, heuristically, we should expect that the larger $\dim X$ is, the smaller $\dim E_s(X)$ should be (i.e. the larger the object, the fewer directions with small shadow). To this end, though it does not have a matching discrete analogue above, it is worth noting Falconer's exceptional set estimate from 1982 \cite{Falconer82}: 
\[
\dim E_s(X) \leq \max\{1 + s - \dim X, 0\}.
\]

Both Kaufman and Falconer's projection theorems were quite difficult to improve upon for general $s$. The fundamental problem in improving upon Kaufman's exceptional set estimate is demonstrated by a hypothetical subring example. If $A \subset \R$ were an $s$-dimensional subring, then for $a\in A$ and $\theta$ in direction $(1, a)$ we have $$P_\theta(A \times A) = c_a(A + aA) = c_a A$$ where $c_a = (1+a^2)^{-1/2}$. Thus, each such $\theta$ is in $E_s(A\times A)$. In other words, $s\leq \dim E_s(A \times A)$, which shows \eqref{eqn:Kaufman} is sharp for this hypothetical example.

However, Bourgain's discretized sum-product theorem \cite{BourgainProjectionTheorem} gives that there does not exist an $s$-dimensional Borel subring of $\R$ for any $s\in (0,1)$. Using this, Orponen and Shmerkin \cite{OrponenShmerkinEpsilonImprove} proved a groundbreaking $\epsilon$-improvement on Kaufman's theorem when $s< \dim X\leq 1$. They proved that for such ranges of $s$ and $t = \dim X$, there exists an $\epsilon:= \epsilon(s,t)$ such that 
\[
\dim E_s(X) \leq s - \epsilon.
\]
This small improvement opened the floodgates on a number of problems in geometric measure theory, ultimately leading to the resolution of Oberlin's conjecture \eqref{eqn:OSRW} due to Orponen--Shmerkin \cite{OrponenShmerkinABC} and Ren--Wang \cite{RenWang} in 2023.

\section*{A Peek Beyond the Floodgates}

Let us conclude by noting a few directions one can go from here. Firstly, it is worth observing that at its core, exceptional set estimates grapple with the question: How can points lie on lines? To this end, exceptional sets are closely related to \textit{Furstenberg sets}. The discrete analogue of the Furstenberg set problem asks: 
Given integers $s\geq 2$ and $t\geq 1$, how small can a finite set $X\subset \R^2$ be if there exists a family $
\mathcal{L}$ of at least $t$ distinct lines such that $|X\cap \ell|\geq s$ for every $\ell\in \mathcal{L}$. 
Using the Szemer\'edi--Trotter theorem, one can prove
\[
|X| \gtrsim \min\{st, s^{3/2}t^{1/2}\}.
\]
The Furstenberg set problem correspondingly asks: 
For $0< s\le 1$ and $0<t\le 2$, given $X\subset \R^2$
such that there exists a set of lines $\mathcal L$ with $\dim \mathcal L= t$ and $\dim (X\cap \ell)\geq s$ for all $\ell \in \mathcal L$, how small can $\dim X$ be? 
Building on work of Orponen--Shmerkin \cite{OrponenShmerkinABC}, Ren--Wang \cite{RenWang} proved the sharp bound for this problem as well:
\[
\dim X \geq \min\{s + t, \frac{3s + t}{2}, s + 1\}.
\]
The Furstenberg set problem and its connection to exceptional set estimates is a rich story that we do not go into further here.

Secondly, it is worth noting that projection theorists study exceptional sets estimates for many types of projections. One important example to this end involves \textit{radial projections} $\pi_x: \R^2 \setminus \{x\}\to \mathbb{S}^1/\pm$ given by
\[
\pi_x(y) = \frac{y-x}{|y-x|}.
\]
Notice that the fibers of these projections are yet again lines, which might lead one to ask if our incidence theoretic approach carries over to this setting. While this line of questioning does give a number of interesting estimates, for more we refer the reader to \cite{OrponenShmerkinABC}.

Finally, in higher dimensions one can consider orthogonal projections onto subspaces of a fixed dimension.  To this end, the reader can verify that \eqref{eqn:CSineq} and \eqref{eqn:STSharp} hold replacing projections onto lines in $\R^2$ with projections onto hyperplanes in $\R^n$, as the fibers of such projections are lines, and thus amenable to the incidence estimates above. In this same vein, Kaufman's estimate \eqref{eqn:Kaufman} holds for hyperplane projections as well, and the sharp analogue of the Szemer\'edi--Trotter-type estimate was recently obtained by Cholak et al. \cite{cholak2026} using \cite{OrponenShmerkinABC,RenWang}. However, the search for sharp exceptional set estimates onto high codimensional planes is still underway, making this an exceptionally exciting time to study projection theory.

\printbibliography

\end{document}